\documentclass[11pt]{article}

\RequirePackage[OT1]{fontenc}

\RequirePackage[aop,amsthm,amsmath,natbib,noinfoline]{imsart_mod}
\RequirePackage[aop]{imsart_mod}

\RequirePackage[dvips,colorlinks]{hyperref}

\usepackage{amssymb,latexsym}

\startlocaldefs

\newtheorem{thm}{Theorem}[section]
\newtheorem{lemma}[thm]{Lemma}

\numberwithin{equation}{section}

\def\al{\alpha}
\def\be{\beta}
\def\ga{\gamma}
\def\Ga{\Gamma}
\def\de{\delta}
\def\ep{\varepsilon}
\def\ze{\zeta}
\def\La{\Lambda}
\def\si{\sigma}
\def\ph{\varphi}
\def\om{\omega}
\def\Om{\Omega}

\def\RR{\mathbb{R}}

\def\FF{\mathcal{F}}
\def\DD{\mathcal{D}}
\def\UU{\mathcal{U}}
\def\YY{\mathcal{Y}}

\def\St{\tilde{S}}

\def\Zh{\hat{Z}}
\def\xih{\hat{\xi}}

\def\To{\Rightarrow}
\def\eqd{\overset{d}{=}}

\providecommand{\flr}[1]{\lfloor#1\rfloor}

\def\pf{\noindent{\bf Proof.} }
\def\qed{{\hfill $\Box$}}

\endlocaldefs

\begin{document}

\begin{frontmatter}

\title{Asymptotic Behavior of a Generalized TCP Congestion Avoidance Algorithm}
\runtitle{Asymptotics of Congestion Avoidance}


\author{\fnms{Teunis J.} \snm{Ott}
\ead[label=u1,url]{http://www.teunisott.com}}
\address{Ott Associates\\ 31 Mountainview Drive\\
Chester, NJ 07930\\
\printead{u1}}%
\affiliation{Ott Associates, Chester, NJ}%

\and

\author{\fnms{Jason} \snm{Swanson}\thanksref{t2}
\ead[label=u2,url]{www.math.wisc.edu/$\sim$swanson}}
\thankstext{t2}{This work was supported in part by the VIGRE
grant of University of Wisconsin-Madison.}
\address{Mathematics Department\\ University of
Wisconsin-Madison\\ 480 Lincoln Dr.\\ Madison, WI 53706-1388\\
\printead{u2}}%
\affiliation{University of Wisconsin-Madison\\ \bigskip
  {\rm August 18, 2006} \bigskip}

\runauthor{T. J. Ott and J. Swanson}

\begin{abstract}

The Transmission Control Protocol (TCP) is a Transport Protocol
used in the Internet. In \citep{O}, a more general class of
candidate Transport Protocols called ``protocols in the TCP
Paradigm'' is introduced. The long run objective of studying this
class is to find protocols with promising performance
characteristics. This paper studies Markov chain models derived
from protocols in the TCP Paradigm.

Protocols in the TCP Paradigm, as TCP, protect the network from
congestion by reducing the ``Congestion Window'' (the amount of
data allowed to be sent but not yet acknowledged) when there is
packet loss or packet marking, and increasing it when there is no
loss. When loss of different packets are assumed to be independent
events and the probability $p$ of loss is assumed to be constant,
the protocol gives rise to a Markov chain $\{W_n\}$, where $W_n$
is the size of the congestion window after the transmission of the
$n$-th packet.

For a wide class of such Markov chains, we prove weak convergence
results, after appropriate rescaling of time and space, as
$p\to0$. The limiting processes are defined by stochastic
differential equations. Depending on certain parameter values, the
stochastic differential equation can define an Ornstein-Uhlenbeck
process or can be driven by a Poisson process.

\end{abstract}



\end{frontmatter}

\section{Introduction}

The Congestion Avoidance algorithm of TCP is designed to prevent
network congestion during the transmission of data over a computer
network. It does this by controlling the congestion window, i.e.
the amount of data ``transmitted but not yet acknowledged'' by a
sender. What follows is a simplified description of a more general
class of Transport Protocols.

Under appropriate units, the congestion window $W$ determines the
maximum amount of data that a source can send without
acknowledgement. The ``TCP Paradigm'' (see \citep{O}) is a class
of protocols that includes TCP (and other Transport Protocols).
For each protocol in the TCP Paradigm there are two functions,
$incr(\cdot)$ and $decr(\cdot)$. If, while the congestion window
equals $W$, a packet is found to be lost (or marked, under ECN --
see \citep{F} and \citep{R}), then the congestion window is
reduced by $decr(W)$. However, the congestion window is never
reduced below some fixed minimum value $\ell\ge0$. If there are no
lost packets, then the congestion window is increased by
$incr(W)$. For protocols in the TCP Paradigm, $incr(W)=c_1 W^\al$
and $decr(W)=c_2 W^\be$. In the special case of TCP, we have
$c_1=1$, $\al=-1$, $c_2=1/2$, and $\be=1$. Another special case of
interest is when $\al=0$ and $\be=1$. This is the algorithm which
Tom Kelly calls ``Scalable TCP'' in \citep{K} and \citep{K2}.

Let $W_n$ denote the size of the congestion window after the
transmission of the $n$-th packet, and let $\chi_n$ be the
indicator function of the event that the $n$-th packet is lost. We
shall assume that the $\chi_n$'s are independent and identically
distributed. In particular, we are assuming that $p=P(\chi_n=1)$
is a constant that does not change with time. Under these
assumptions, we are led to the parameterized family of Markov
processes
  \begin{equation}\label{Wn}
  W_{p,n+1} = (W_{p,n} + c_1W_{p,n}^\al(1 - \chi_{p,n+1})
    - c_2W_{p,n}^\be\chi_{p,n+1}) \vee \ell.
  \end{equation}
The assumptions we place on the various parameters in the model
are:
  \begin{align}
  & \{\chi_{p,n}\}_{n=1}^\infty\text{ is an iid sequence of }
    \{0,1\}\text{-valued random variables,}\\
  & p = P(\chi_{p,n} = 1),\\
  & c_1 > 0 \text{ and } c_2 > 0,\\
  & -\infty < \al < \be \le 1 \text{ and } \ell \ge 0,\\
  & \text{if } \be = 1, \text{ then } c_2 < 1, \text{ and}\\
  & \text{if } \be < 1, \text{ then } \ell > 0.
  \end{align}
We will frequently drop the dependence on $p$ from our notation
and simply refer to the processes $\{\chi_n\}$ and $\{W_n\}$.

We are interested in studying the asymptotic behavior of $\{W_n\}$
as $p\to0$. To this end, we define the continuous time process
  \begin{equation}\label{Zp}
  Z_p(t) = p^\ga W_{\flr{tp^{-\nu}}}, \text{ where }
    \ga = (\be - \al)^{-1} \text{ and } \nu = (1 - \al)\ga.
  \end{equation}
In the case that $\be=1$, we will show that $Z_p$ converges weakly
as $p\to0$ to the process $Z$ defined by
  \begin{equation}\label{Z}
  Z(t) = Z(0) + c_1\int_0^t Z(s)^\al\,ds
    - c_2\int_0^t Z(s-)\,dN(s),
  \end{equation}
where $N$ is a unit rate Poisson process, independent of
$Z(0)=\lim Z_p(0)$. (Note that this is the conjecture given on
page 362 of \citep{O}.) We will also show that, when $\ell>0$, the
stationary distributions of the discrete time Markov chains
$\{p^\ga W_n\}$ converge weakly to the unique stationary
distribution of $Z$. Questions about the convergence of the
stationary distributions when $\be=1$, as well as the rate of
convergence, are addressed in \citep{O2} using techniques that
differ from those used in this paper.

In the case that $\be<1$, we will show that $Z_p$ converges to the
process $\ze$ defined by
  \begin{equation}\label{ze}
  \ze(t) = \ze(0) + \int_0^t
    (c_1\ze(s)^\al - c_2\ze(s)^\be)\,ds,
  \end{equation}
where $\ze(0)=\lim Z_p(0)$. With the exception of the initial
condition, the process $\ze$ is entirely deterministic. The
convergence of $Z_p$ to $\ze$ is therefore a law of large numbers
type of result. Hence, in the case $\be<1$, we can extend our
analysis and study the fluctuations of $Z_p$ around this central
tendency. Unfortunately, it will not suffice to center $Z_p$ by
$\ze$. We must rather define
  \begin{equation}\label{zep}
  \ze_p(t) = \ze_p(0) + \int_0^t
    (c_1(1 - p)\ze_p(s)^\al - c_2\ze_p(s)^\be)\,ds,
  \end{equation}
where $\ze_p(0)\to\ze(0)$, and consider the processes
  \begin{equation}\label{xip}
  \xi_p(t) = p^{-\tau}(Z_p(t) - \ze_p(t)), \text{ where }
    \tau = (\nu - 1)/2.
  \end{equation}
We will show that $\xi_p$ converges weakly as $p\to0$ to the
process $\xi$ defined by
  \begin{equation}\label{xi}
  \begin{split}
  \xi(t) &= \xi(0)
    + \int_0^t (c_1\al\ze(s)^{\al-1}
    - c_2\be\ze(s)^{\be-1})\xi(s)\,ds\\
    &\quad - c_2\int_0^t \ze(s)^\be\,dB(s),
  \end{split}
  \end{equation}
where $B$ is a Brownian motion and $\xi(0)=\lim\xi_p(0)$.

A special case of this last result is worth mentioning. For each
$p\in[0,1)$, define
  \begin{equation}\label{cp}
  c_p = (c_1(1-p)/c_2)^\ga,
  \end{equation}
so that $\ze_p(t)=c_p$ is an invariant solution to \eqref{zep}.
Also, $\ze(0)=\lim\ze_p(0)=c_0$ is an invariant solution to
\eqref{ze}. Hence, for an appropriate choice of $Z_p(0)$, $\xi_p$
converges to the Ornstein-Uhlenbeck process defined by
  \begin{equation}\label{OU}
  d\xi = -\mu\xi dt + \si dW,
  \end{equation}
where $W=-B$,
  \begin{align*}
  \mu &= c_2\be(c_1/c_2)^{\ga(\be-1)}
    - c_1\al(c_1/c_2)^{\ga(\al-1)}\\
  &= (\be - \al)c_1^{-(1-\be)/(\be-\al)}c_2^{(1-\al)/(\be-\al)},
  \end{align*}
and
  \[
  \si = c_2(c_1/c_2)^{\ga\be}
    = c_1^{\be/(\be-\al)}c_2^{-\al/(\be-\al)}.
  \]
(Note that this is the conjecture given on page 364 of \citep{O}.)
We will also show that the stationary distributions of the
discrete time Markov chains $\{p^{-\tau}(p^\ga W_n-c_p)\}$
converge weakly to the unique stationary distribution of the above
Ornstein-Uhlenbeck process.

\section{Main Results}

We first consider the case $\be=1$ and begin by cataloging some
properties of the limit process $Z$.

\begin{lemma}\label{L:Zprop}
If $Z(0)>0$ a.s., then the stochastic differential equation
\eqref{Z} has a unique solution $Z$. With probability one,
$Z(t)>0$ for all $t\ge0$. Moreover, if
$\tau=\inf\{t\ge0:Z(t)=c_0\}$, where $c_0$ is given by \eqref{cp},
then $\tau<\infty$ a.s.
\end{lemma}

\pf For each realization of the Poisson process, \eqref{Z} can be
solved deterministically and the solution is unique. Let
  \[
  T = \inf\{t \ge 0: Z(t) \notin (0,\infty)\}.
  \]
Since $Z$ decreases only at the jump times of the Poisson process,
and, with probability one, these jump times have no accumulation
points, it follows that $T=\infty$ a.s.

To show that $\tau<\infty$ a.s., it will suffice to assume that
$Z(0)=x>0$ is deterministic. We first consider the case $x\le
c_0$. Suppose $\tau(\om)=\infty$. Then $Z(t,\om)<c_0$ for all
$t\ge0$. Find $u>r$ such that $u-r>\ga c_2^{-1}$ and
$N(u,\om)=N(r,\om)$. Then for all $t\in(r,u]$,
  \[
  Z(t,\om) = Z(r,\om) + c_1\int_r^t Z(s,\om)^\al\,ds.
  \]
Since the solution to this integral equation is unique,
  \[
  Z(t,\om) = (c_1(1 - \al)(t - r) + Z(r,\om)^{1-\al})^\ga.
  \]
Therefore,
  \[
  c_0 > Z(u,\om) > (c_1(1 - \al)(u - r))^\ga > c_0,
  \]
a contradiction. Hence, $\tau<\infty$ a.s.

We next consider the case $x>c_0$. Define
  \[
  \si_1 = \inf\{t \ge 0 : Z(t) < c_0\} \text{ and }
  \si_2 = \inf\{t \ge \si_1 : Z(t) = c_0\},
  \]
so that $\tau\le\si_2$, and it will suffice to show that
$\si_2<\infty$ a.s. Fix $L>x$ and define
$\rho=\inf\{t\ge0:Z(t)\notin[c_0,L]\}$. Suppose
$\rho(\om)=\infty$. Then $Z(t,\om)\in[c_0,L]$ for all $t\ge0$. Let
  \[
  K = \inf\{u^\al : c_0 \le u \le L\} > 0.
  \]
Find $u>r$ such that $u-r>(L-c_0)/(c_1K)$ and $N(u,\om)=N(r,\om)$.
Then
  \[
  L \ge Z(u,\om) = Z(r,\om)
    + c_1\int_r^u Z(s,\om)^\al\,ds \ge c_0 + c_1(u - r)K > L,
  \]
a contradiction. Hence, $\rho<\infty$ a.s.

Now, observe that
  \[
  Z(t\wedge\rho) = x
    + \int_0^{t\wedge\rho} (c_1Z(s)^\al - c_2Z(s))\,ds
    - c_2\int_0^{t\wedge\rho} Z(s-)\,dM(s),
  \]
where $M(t)=N(t)-t$ is the compensated Poisson process. If
$s<t\wedge\rho$, then $Z(s)\ge c_0=(c_1/c_2)^\ga$. This implies
that $c_1Z(s)^\al-c_2Z(s)\le0$. Since $M$ is a martingale,
$E[Z(t\wedge\rho)]\le x$. Letting $t\to\infty$ gives
$E[Z(\rho)]\le x$. Hence, $P(Z(\rho)=L)\le x/L$. Note that either
$Z(\rho)=L$ or $Z(\rho)<c_0$. Therefore,
  \[
  P(\si_1 = \infty) \le P(Z(\rho) = L) \le x/L.
  \]
Letting $L\to\infty$ shows $\si_1<\infty$ a.s.

As in Theorem V.6.35 in \citep{P}, $Z$ is a strong Markov process.
Therefore,
  \[
  P(\si_2 = \infty) = E[P^{Z(\si_1)}(\tau = \infty)].
  \]
But $Z(\si_1)<c_0$, and we have already shown that
$P^x(\tau=\infty)=0$ for all $x\le c_0$. Hence, $\si_2<\infty$
a.s. \qed

We are now prepared to state our main results for the case
$\be=1$. If $\mu_p$ and $\mu$ are Borel measures on a metric space
$S$, then the notation $\mu_p\To\mu$ will mean that $\mu_p$
converges weakly to $\mu$ as $p\to0$, that is, $\int_S
f\,d\mu_p\to\int_S f\,d\mu$ as $p\to0$ for all bounded, continuous
$f:S\to\RR$. If $X_p$ and $X$ are $S$-valued random variables,
then $X_p\To X$ will mean that $PX_p^{-1}\To PX^{-1}$. When $X_p$
and $X$ are processes, we will take our metric space to be
$D_{\RR^d}[0,\infty)$, the space of cadlag functions from
$[0,\infty)$ to $\RR^d$, with the Skorohod metric. See \citep{EK}
for details.

\begin{thm}\label{T:Poisson}
Suppose $\be=1$. Let the processes $Z_p$ be given by \eqref{Zp}
and suppose that $Z_p(0)\To Z(0)$, where $Z(0)>0$ a.s. Let $Z$ be
the unique solution to \eqref{Z}. Then $Z_p\To Z$.
\end{thm}

\begin{thm}\label{T:be=1}
Suppose $\be=1$ and $\ell>0$. Then the Markov chain $\{W_n\}$ has
a unique stationary distribution. Moreover, the process $Z$ given
by \eqref{Z} has a unique stationary distribution $\eta$ on
$(0,\infty)$. For each $p>0$, let $\eta_p$ be the stationary
distribution for the Markov chain $\{p^\ga W_n\}$. Then
$\eta_p\To\eta$.
\end{thm}

For some results on stationary distributions in the case $\be=1$
and $\ell=0$, see \citep{O3}.

For the case $\be<1$, we need some preliminary definitions. Assume
that for all $p\in(0,1)$, the processes $\{W_{p,n}\}$ are defined
on the same probability space $(\Om,\FF,P)$. Define the
$\si$-algebra
  \begin{equation}\label{F0}
  \FF_0 = \si(W_{p,0}:0<p<1) \vee \mathcal{N},
  \end{equation}
where $\mathcal{N}$ denotes the collection of events $D\in\FF$
with $P(D)=0$.

\begin{thm}\label{T:LLN}
Suppose $\be<1$. Let the processes $Z_p$ be given by \eqref{Zp}.
Suppose that $Z_p(0)\To\ze(0)$, where $\ze(0)>0$ a.s. Let $\ze$
the unique solution to \eqref{ze}. Then $Z_p\To\ze$. Moreover, if
$Z_p(0)\to\ze(0)$ in probability, then $Z_p\to\ze$ in probability.
\end{thm}

\begin{thm}\label{T:CLT}
Suppose $\be<1$. Let the processes $Z_p$ be given by \eqref{Zp}.
For each $p\in(0,1)$, let $\ze_p(0)$ be a strictly positive random
variable defined on $(\Om,\FF,P)$. Assume that $\ze_p(0)$ is
$\FF_0$-measurable and $Z_p(0)-\ze_p(0)\to0$ in probability.
Define $\ze_p$ and $\xi_p$ by \eqref{zep} and \eqref{xip},
respectively.

Suppose that there exists a pair of random variables
$(\xi(0),\ze(0))$, defined on $(\Om,\FF,P)$, such that $\ze(0)>0$
a.s., $\ze_p(0)\to\ze(0)$ in probability, and
$(\xi_p(0),\ze_p(0))\To(\xi(0),\ze(0))$. Let $B$ be a standard
Brownian motion independent of $(\xi(0),\ze(0))$ and define the
processes $\ze$ and $\xi$ by \eqref{ze} and \eqref{xi},
respectively. Then $(\xi_p,\ze_p)\To(\xi,\ze)$.
\end{thm}

\begin{thm}\label{T:be<1}
Suppose $\be<1$. Then the Markov chain $\{W_n\}$ has a unique
stationary distribution. For each $p>0$, let $\eta_p$ be the
stationary distribution for the Markov chain $\{p^{-\tau}(p^\ga
W_n-c_p)\}$. Then $\eta_p\To\eta$, where $\eta$ is the stationary
distribution of the Ornstein-Uhlenbeck process given by
\eqref{OU}.
\end{thm}

\section{General Definitions}

Define
  \[
  \La_n = (\ell - W_{n-1} - c_1W_{n-1}^\al(1 - \chi_n)
    + c_2W_{n-1}^\be\chi_n) \vee 0,
  \]
so that
  \[
  W_{n+1} = W_n + c_1 W_n^\al
    - (c_1 W_n^\al + c_2 W_n^\be)\chi_{n+1} + \La_{n+1}.
  \]
If we let $W(t)=W_{\flr{t}}$, then we can rewrite this recursive
relation as the integral equation
  \begin{equation*}
  \begin{split}
  W(t) &= W(0) + c_1\int_0^t W(s-)^\al\,dm(s)\\
  &\quad - \int_0^t (c_1 W(s-)^\al + c_2 W(s-)^\be)\,dS(s) +
  L(t),
  \end{split}
  \end{equation*}
where
  \[
  m(t) = \flr{t},\,\, S(t) = \sum_{j=1}^{\flr{t}}\chi_j,
    \text{ and } L(t) = \sum_{j=1}^{\flr{t}}\La_j.
  \]
Using \eqref{Zp}, it is then easy to see that
  \begin{equation}\label{ZpIE}
  \begin{split}
  Z_p(t) &= Z_p(0) + c_1\int_0^t Z_p(s-)^\al\,dm_p(s)\\
  &\quad - c_1p\int_0^t Z_p(s-)^\al\,dS_p(s)
    - c_2\int_0^t Z_p(s-)^\be\,dS_p(s) + L_p(t),
  \end{split}
  \end{equation}
where
  \[
  m_p(t) = p^\nu m(tp^{-\nu}),\,\,
  S_p(t) = p^{\nu-1}S(tp^{-\nu}), \text{ and }
  L_p(t) = p^\ga L(tp^{-\nu}).
  \]
Note that if we define the filtration
  \[
  \FF_t^p = \FF_0 \vee \si(\chi_{p,j} : j \le \flr{tp^{-\nu}}),
  \]
then $m_p$, $S_p$, and $L_p$ are all $\{\FF_t^p\}$-adapted.

Define the $\RR^2$-valued cadlag $\{\FF_t^p\}$-semimartingale
  \[
  Y_p = (m_p,S_p)^T
  \]
and define the function $G_p:\RR^2\to\RR$ by
  \[
  G_p(x) = (c_1 x^\al, -c_1 p x^\al - c_2 x^\be)1_{\{x>0\}}.
  \]
Then \eqref{ZpIE} becomes
  \[
  Z_p(t) = Z_p(0) + \int_0^t G_p(Z_p(s-))\,dY_p(s) + L_p(t).
  \]
To show that $Z_p$ converges as $p\to0$, we will apply the
theorems in \citep{KP}. This approach, however, comes with two
technical difficulties. The first is the presence of the local
time term $L_p$; the second is the fact that $G_p$ may have a
singularity at the origin. To deal with these issues, we introduce
the process $Z_p^\ep$, defined as the unique solution to
  \begin{equation}\label{Zpep}
  Z_p^\ep(t) = Z_p(0) + \int_0^t G_p^\ep(Z_p^\ep(s-))\,dY_p(s),
  \end{equation}
where $G_p^\ep=G_p(\ep)1_{(-\infty,\ep)}+G_p1_{[\ep,\infty)}$. To
quantify the sense in which $Z_p$ and $Z_p^\ep$ are close, we
define the functional $h_\ep:D_{\RR^d}[0,\infty)\to[0,\infty]$ by
  \[
  h_\ep(x) = \inf\{t \ge 0: |x(t)| \wedge |x(t-)| \le \ep\},
  \]
and the stopping times $\tau_p(\ep)=h_\ep(Z_p^\ep)$, and we
observe that
  \begin{equation}\label{local}
  L_p = 0 \text{ and } Z_p = Z_p^\ep \text{ on }
    [0,\tau_p(\ep\vee p^\ga\ell)).
  \end{equation}
By (3.5.2) in \citep{EK}, if two cadlag functions $x$ and $y$
agree on the interval $[0,t)$, then $d(x,y)\le e^{-t}$, where $d$
is the metric on $D_{\RR^d}[0,\infty)$.

\section{Convergence of $Z_p$}

In this section, we will prove Theorems \ref{T:Poisson} and
\ref{T:LLN} by applying the theorems in \citep{KP} to the
processes $Z_p^\ep$ given by \eqref{Zpep}. We must therefore
define the processes to which they converge in the cases $\be=1$
and $\be<1$.

Let $G(x)=(c_1 x^\al,-c_2x^\be)1_{\{x>0\}}$ and
$G^\ep=G(\ep)1_{(-\infty,\ep)}+G1_{[\ep,\infty)}$, and note that
$G_p^\ep\to G^\ep$ uniformly on compacts as $p\to0$. Let $N$ be a
unit rate Poisson process, define
  \[
  Y(t) = (t,N(t))^T \text{ and } y(t) = (t,t)^T,
  \]
and let $Z^\ep$ and $\ze^\ep$ be the unique solutions to
  \begin{align}
  Z^\ep(t) &= Z(0)
    + \int_0^t G^\ep(Z^\ep(s-))\,dY(s),\label{Zep}\\
  \ze^\ep(t) &= \ze(0)
    + \int_0^t G^\ep(\ze^\ep(s-))\,dy(s),\label{zeep}
  \end{align}
where $Z(0)$ and $N$ are independent. Note that if $\be=1$, then
$Z^\ep=Z$ on $[0,h_\ep(Z^\ep))$ and
$h_\ep(Z^\ep)=h_\ep(Z)\to\infty$ a.s. as $\ep\to0$. Hence,
$d(Z^\ep,Z)\le\exp(-h_\ep(Z))\to0$ a.s. That is, $Z^\ep\to Z$ a.s.
in $D_\RR[0,\infty)$. Similarly, if $\be<1$, then $\ze^\ep=\ze$ on
$[0,h_\ep(\ze^\ep))$, $h_\ep(\ze^\ep)=h_\ep(\ze)\to\infty$ a.s.,
and $\ze^\ep\to\ze$ a.s. in $D_\RR[0,\infty)$.

We will show that $Z_p^\ep\To Z^\ep$ and $\ze_p^\ep\To\ze^\ep$. To
pass from this to the conclusions of Theorems \ref{T:Poisson} and
\ref{T:LLN}, we will need the following lemma, which is easily
proved using the Prohorov metric. (See Section 3.1 in \citep{EK}.)

\begin{lemma}\label{L:Proho}
Let $(S,d)$ be a complete and separable metric space. Let
$\{X_p\}_{p>0}$ be a family of $S$-valued random variables and
suppose, for each $\ep$, there exists a family $\{X_p^\ep\}_{p>0}$
such that
  \[
  \limsup_{p\to0}E[d(X_p,X_p^\ep)] \le \de_\ep,
  \]
where $\de_\ep\to0$ as $\ep\to0$. Suppose also that for each
$\ep$, there exists $Y^\ep$ such that $X_p^\ep\To Y^\ep$ as
$p\to0$. Then there exists $X$ such that $X_p\To X$ and $Y^\ep\To
X$.
\end{lemma}

\noindent{\bf Proof of Theorem \ref{T:Poisson}.} Suppose $\be=1$,
$Z_p$ is given by \eqref{Zp}, and $Z_p(0)\To Z(0)$, where $Z(0)>0$
a.s. Let $Z$ be the solution to \eqref{Z}.

Let $Z_p^\ep$ and $Z^\ep$ be given by \eqref{Zpep} and
\eqref{Zep}. We first show that $Z_p^\ep\To Z^\ep$. Recall that
$G_p^\ep\to G^\ep$ uniformly on compacts. Also observe that
$S_p\To N$ (see, for example, Problem 7.1 in \citep{EK}). Hence,
since $Z_p(0)$ and $Y_p$ are independent,
$(Z_p(0),Y_p)\To(Z(0),Y)$ in $D_{\RR^3}[0,\infty)$. Hence, by
Theorem 5.4 in \citep{KP}, it will suffice to show that $Y_p$ has
a semimartingale decomposition $Y_p=M_p+A_p$ into a martingale
part and a bounded variation part such that for each $t\ge0$,
  \begin{equation}\label{KPcond}
  \sup_p E[[M_p]_t + T_t(A_p)] < \infty,
  \end{equation}
where $[M_p]_t$ is the quadratic variation process of $M_p$ and
$T_t(A_p)$ is the total variation of $A_p$ on the interval
$[0,t]$.

For this, define
  \[
  \St_p(t) = S_p(t) - m_p(t)
    = p^{\nu-1}\sum_{j=1}^{\flr{tp^{-\nu}}}(\chi_j - p),
  \]
so that $\St_p$ is an $\{\FF_t^p\}$-martingale. Note that
$T_t(m_p)=m_p(t)$ and
  \begin{equation}\label{qv}
  E[\St_p]_t
    = p^{2\nu-2}\sum_{j=1}^{\flr{tp^{-\nu}}}E|\chi_j - p|^2
    = p^{2\nu-2}\flr{tp^{-\nu}}p(1 - p)
    \le tp^{\nu-1}.
  \end{equation}
Since $\be=1$ implies $\nu=1$, this verifies \eqref{KPcond} and
shows that $Z_p^\ep\To Z^\ep$.

By passing to a subsequence, we can assume there exists a
$[0,\infty]$-valued random variable $\si(\ep)$ such that
$(Z_p^\ep,h_\ep(Z_p^\ep))\To(Z^\ep,\si(\ep))$. By \eqref{local},
  \begin{align*}
  \limsup_{p\to0} E[d(Z_p,Z_p^\ep)]
    &\le \limsup_{p\to0} E[\exp(-\tau_p(\ep\vee p^\ga\ell))]\\
  &= \limsup_{p\to0} E[\exp(-h_\ep(Z_p^\ep))]\\
  &= E[\exp(-\si(\ep))].
  \end{align*}
We claim that $E[\exp(-\si(\ep))]\le E[\exp(-h_\ep(Z^\ep))]$. To
see this, let us assume by the Skorohod Representation Theorem
(see Theorem 3.1.8 in \citep{EK}) that
$(Z_p^\ep,h_\ep(Z_p^\ep))\to(Z^\ep,\si(\ep))$ a.s. Then
$h_\ep(Z^\ep)\le\si(\ep)$ a.s., which proves the claim.

Since $h_\ep(Z^\ep)=h_\ep(Z)\to\infty$ a.s. as $\ep\to0$, we can
apply Lemma \ref{L:Proho} to conclude that $Z_p\To Z$. \qed

\bigskip

\noindent{\bf Proof of Theorem \ref{T:LLN}.} Suppose $\be<1$,
$Z_p$ is given by \eqref{Zp}, and $Z_p(0)\To\ze(0)$, where
$\ze(0)>0$ a.s. Let $\ze$ be the solution to \eqref{ze}.

Note that $\be<1$ implies $\nu>1$. Hence, \eqref{qv} implies that
\eqref{KPcond} is satisfied and $\St_p\to0$ in probability.
Therefore, $(Z_p(0),Y_p)\To(Z(0),y)$ in $D_{\RR^3}[0,\infty)$. By
Theorem 5.4 in \citep{KP}, $Z_p^\ep\To\ze^\ep$. By Corollary 5.6
in \citep{KP}, if $Z_p(0)\to\ze(0)$ in probability, then
$Z_p^\ep\to\ze^\ep$ in probability. By the same argument as above,
this implies that $Z_p$ converges to $\ze$ in distribution or in
probability, respectively. \qed

\section{Fluctuations of $Z_p$}

In this section, we prove Theorem \ref{T:CLT}. Let us first recall
the setting of that theorem. We have $\be<1$ and $Z_p$ given by
\eqref{Zp}. Recall that the processes $Z_p$ are all defined on the
same probability space $(\Om,\FF,P)$. For each $p>0$, $\ze_p(0)$
is an $\FF_0$-measurable random variable, where $\FF_0$ is given
by \eqref{F0}, such that $\ze_p(0)>0$ a.s. and
$Z_p(0)-\ze_p(0)\to0$ in probability. The processes $\ze_p$ and
$\xi_p$ are then given by \eqref{zep} and \eqref{xip}.

To apply the theorems in \citep{KP}, we wish to write $\xi_p$ as
the solution to a stochastic differential equation. By \eqref{zep}
and \eqref{ZpIE}, we have
  \begin{equation}\label{xipIE}
  \begin{split}
  \xi_p(t) &= \xi_p(0) + c_1(1 - p)\int_0^t
    p^{-\tau}(Z_p(s-)^\al - \ze_p(s)^\al)\,dm_p(s)\\
  &\quad
    - c_2\int_0^t p^{-\tau}(Z_p(s-)^\be - \ze_p(s)^\be)\,dS_p(s)\\
  &\quad - c_2\int_0^t \ze_p(s)^\be\,dB_p(s) + R_p(t),
  \end{split}
  \end{equation}
where
  \[
  B_p(t) = p^{-\tau}(S_p(t) - m_p(t))
    = p^{(\nu-1)/2}\sum_{j=1}^{\flr{tp^{-\nu}}}(\chi_j - p)
  \]
and
  \begin{equation}\label{Rp}
  \begin{split}
  R_p(t) &= p^{-\tau}\int_0^t
    (c_1(1 - p)\ze_p(s)^\al - c_2\ze_p(s)^\be)\,d(m_p(s) - s)\\
  &\quad - c_1p\int_0^t Z_p(s-)^\al\,dB_p(s) + p^{-\tau}L_p(t).
  \end{split}
  \end{equation}
Given a real number $r$, let us define the continuous function
$F_r:(0,\infty)^2\to\RR$ by
  \[
  F_r(x,y) = \frac{x^r - y^r}{x - y}1_{\{x\ne y\}}
    + ry^{r-1}1_{\{x=y\}}.
  \]
Using this, \eqref{xipIE} becomes
  \begin{equation}\label{xipIE2}
  \begin{split}
  \xi_p(t) &= \xi_p(0) + c_1(1 - p)\int_0^t
    \xi_p(s-)\DD_p^\al(s-)\,dm_p(s)\\
  &\quad - c_2\int_0^t \xi_p(s-)\DD_p^\be(s-)\,dS_p(s)
    - c_2\int_0^t \ze_p(s)^\be\,dB_p(s) + R_p(t),
  \end{split}
  \end{equation}
where $\DD_p^r=F_r(Z_p,\ze_p)$.

\bigskip

\noindent{\bf Proof of Theorem \ref{T:CLT}.} Suppose that there
exists a pair of random variables $(\xi(0),\ze(0))$, defined on
$(\Om,\FF,P)$, such that $\ze(0)>0$ a.s., $\ze_p(0)\to\ze(0)$ in
probability, and $(\xi_p(0),\ze_p(0))\To(\xi(0),\ze(0))$. Since
the map that takes a point $x>0$ to the unique solution to
\eqref{zep} with $\ze_p(0)=x$ is continuous, $\ze_p\to\ze$ in
probability and $(\xi_p(0),\ze_p)\To(\xi(0),\ze)$. Also, since
$F_r$ is continuous, $\DD_p^r\to r\ze(\cdot)^{r-1}$ in
probability.

Let
  \begin{align*}
  \UU_p(t) &= \xi_p(0) - c_2\int_0^t \ze_p(s)^\be\,dB_p(s)
    + R_p(t), \text{ and}\\
  \YY_p(t) &= c_1(1 - p)\int_0^t \DD_p^\al(s-)\,dm_p(s)
    - c_2\int_0^t \DD_p^\be(s-)\,dS_p(s),
  \end{align*}
so that \eqref{xipIE2} becomes
  \begin{equation}\label{xipIE3}
  \xi_p(t) = \UU_p(t) + \int_0^t \xi_p(s-)\,d\YY_p(s).
  \end{equation}
We will apply the theorems in \citep{KP} to this integral
equation.

We first show that $R_p\to0$ in probability. By the Martingale
Central Limit Theorem (Theorem 7.1.4 in \citep{EK}), $B_p\To B$,
where $B$ is a standard Brownian motion; by Theorem \ref{T:LLN},
$Z_p\to\ze$ in probability; and by \eqref{qv}, $\{B_p\}$ satisfies
\eqref{KPcond}. Hence, by Theorem 2.2 in \citep{KP},
  \[
  c_1p\int_0^t Z_p(s-)^\al\,dB_p(s) \to 0
  \]
in probability. By \eqref{local}, $p^{-\tau}L_p=0$ on
$[0,h_{p^\ga\ell}(Z_p))$. Since $h_{p^\ga\ell}(Z_p)\to\infty$ in
probability, $p^{-\tau}L_p\to0$ in probability.

For the final term in \eqref{Rp}, note that
$p^{-\tau}|m_p(t)-t|\le p^{\nu-\tau}$ and $\nu-\tau=(\nu+1)/2>0$.
Hence, $p^{-\tau}(m_p(t)-t)\to0$ uniformly. Let
$f_p(s)=c_1(1-p)\ze_p(s)^\al-c_2\ze_p(s)^\be$. Since $\ze_p\to\ze$
in probability, we can pass to a subsequence and assume that
$\ze_p\to\ze$ uniformly on $[0,t]$, a.s. By \eqref{zep}, this
implies that $\ze_p'\to\ze'$ uniformly on $[0,t]$. Hence, $f_p$
and $f_p'$ converge uniformly. Integrating by parts, we have
  \begin{equation*}
  \begin{split}
  p^{-\tau}\int_0^t f_p(s)\,d(m_p(s) - s)
    &= p^{-\tau}f_p(t)(m_p(t) - t)\\
  &\quad - p^{-\tau}\int_0^t (m_p(s) - s)f_p'(s)\,ds,
  \end{split}
  \end{equation*}
which goes to zero uniformly and completes the proof that
$R_p\to0$ in probability.

It now follows from Theorem 5.2 in \citep{KP} that
$(\UU_p,\YY_p,\ze_p)\To(\UU,\YY,\ze)$, where
  \begin{align*}
  \UU(t) &= \xi(0) - c_2\int_0^t \ze(s)^\be\,dB(s), \text{ and}\\
  \YY(t) &= c_1\int_0^t \al\ze(s)^{\al-1}\,ds
    - c_2\int_0^t \be\ze(s)^{\be-1}\,ds,
  \end{align*}
and $B$ is a standard Brownian motion independent of
$(\xi(0),\ze(0))$. By Remark 2.5 in \citep{KP}, we may apply
Theorem 5.4 in \citep{KP} to \eqref{xipIE3} and conclude that
$(\xi_p,\ze_p)\To(\xi,\ze)$, where $\xi$ is the unique solution to
\eqref{xi}. \qed

\section{Stationary Distributions}

In this section, we prove Theorems \ref{T:be=1} and \ref{T:be<1}.
For this, we make time continuous in a slightly different manner
than before. Let $N$ be a unit rate Poisson process independent of
$\{W_n\}$ and let $X(t)=W_{N(t)}$. Then $X$ is a continuous time
Markov chain on $E=[\ell,\infty)$ with generator
  \[
  A\ph(x) = p(\ph(x - g(x)) - \ph(x))
    + (1 - p)(\ph(x + c_1x^\al) - \ph(x)),
  \]
where $g(x)=(c_2x^\be)\wedge(x-\ell)$. When $\be=1$, we will study
the process
  \[
  \Zh_p(t) = p^\ga X(tp^{-1}),
  \]
whereas when $\be<1$, we will consider
  \[
  \xih_p(t) = p^{-\tau}(p^\ga X(tp^{-\nu}) - c_p),
  \]
where $c_p$ is given by \eqref{cp}. It is easy to see that a
probability measure is a stationary distribution for $\{p^\ga
W_n\}$ or $\{p^{-\tau}(p^\ga W_n-c_p)\}$ if and only if it is a
stationary distribution for $\Zh_p$ or $\xih_p$, respectively.

\begin{lemma}
If $\ell>0$, then $\{W_n\}$ has a unique stationary distribution.
\end{lemma}

\pf It will suffice to show that $X$ has a unique stationary
distribution. Let $\ph(x)=x$ so that
  \[
  A\ph(x) = -pg(x) + (1 - p)c_1x^\al.
  \]
Since $g(x)=c_2x^\be$ for $x$ sufficiently large, $A\ph$ is
bounded above and $A\ph(x)\to-\infty$ as $x\to\infty$. By Lemmas
4.9.5 and 4.9.7 in \citep{EK}, the family of probability measures
$\{\mu_t\}_{t\ge1}$ defined by
  \[
  \mu_t(\Ga) = \frac1t\int_0^t P^x(X(s) \in \Ga)\,ds
  \]
is relatively compact. By Theorem 4.9.3 in \citep{EK}, any
subsequential weak limit of $\{\mu_t\}$ is a stationary
distribution for $X$.

To show that the stationary distribution is unique, it will
suffice to show that for all $x\in E$,
  \[
  \tau = \inf\{t \ge 0 : X(t) = \ell\} < \infty,
    \quad P^x \text{-a.s.}
  \]
(See, for example, Problem 4.36 in \citep{EK}.) Let $x\in E$ be
arbitrary and let $\ep>0$. Choose $M$ such that
$\mu_t([\ell,M])\ge1-\ep$ for all $t\ge0$. Note that there exists
$K>0$ such that $P^y(\tau<\infty)\ge K$ for all $y\in[\ell,M]$.

Define the stopping times $\tau_0=0$ and
  \[
  \tau_{j+1} = \inf\{t \ge \tau_j+1 : X(t) \le M\},
  \]
and note that $\tau_j\to\infty$ a.s. By the strong Markov
property,
  \begin{align*}
  P(\tau = \infty, \tau_j < \infty)
    &= E[1_{\{\tau \ge \tau_j, \tau_j < \infty\}}
    P^{X(\tau_j)}(\tau = \infty)]\\
  &\le (1 - K)P(\tau \ge \tau_j, \tau_j < \infty)
  \end{align*}
Letting $j\to\infty$ shows that $P(\{\tau=\infty\}\cap D)=0$,
where $D$ is the event that $\tau_j<\infty$ for all $j$. Note that
  \[
  1_{D^c} \le
    \liminf_{t\to\infty}\frac1t
    \int_0^t 1_{\{X(s) > M\}}\,ds.
  \]
Hence, by Fatou's Lemma,
$P(D^c)\le\liminf_{t\to\infty}\mu_t((M,\infty))\le\ep$. Therefore,
$P(\tau=\infty)=P(\{\tau=\infty\}\cap D^c)\le\ep$. Since $\ep$ was
arbitrary, $\tau<\infty$ $P^x$-a.s. and the stationary
distribution is unique. \qed

\bigskip

\noindent{\bf Proof of Theorem \ref{T:be=1}.} In what follows, $C$
and $K$ will denote strictly positive, finite constants that do
not depend on $p$ and may change value from line to line.

Suppose $\be=1$, $\ell>0$, and $\eta_p$ is the stationary
distribution for $\{p^\ga W_n\}$. Then $\eta_p$ is the stationary
distribution for $\Zh_p$, which is a continuous time Markov chain
on $E_p=[p^\ga\ell,\infty)$ with generator
  \begin{equation*}
  \begin{split}
  A_p\ph(x) &= \ph(x - p^\ga g(p^{-\ga}x)) - \ph(x)\\
  &\quad + p^{-1}(1 - p)(\ph(x + pc_1x^\al) - \ph(x)).
  \end{split}
  \end{equation*}
Let $\ph(x)=x+x^{-1}$, so that
  \begin{equation*}
  \begin{split}
  A_p\ph(x) &= -p^\ga g(p^{-\ga}x) + (1 - p)c_1x^\al\\
    &\quad + \frac{p^\ga g(p^{-\ga}x)}{x(x - p^\ga g(p^{-\ga}x))}
    - \frac{(1 - p)c_1x^\al}{x(x + pc_1x^\al)}.
  \end{split}
  \end{equation*}
Since $x\mapsto1+pc_1x^{\al-1}$ is decreasing,
  \[
  1 + pc_1x^{\al-1} \le 1 + pc_1(p^\ga\ell)^{\al-1}
    = 1 + c_1\ell^{\al-1}
  \]
for all $x\in E_p$. Hence,
  \[
  A_p\ph(x) \le -p^\ga g(p^{-\ga}x) + Cx^\al
    + \frac{p^\ga g(p^{-\ga}x)}{x(x - p^\ga g(p^{-\ga}x))}
    - Kx^{\al-2}
  \]
whenever $p<1/2$.

If $x\ge p^\ga\ell/(1-c_2)$, then $g(p^{-\ga}x)=c_2p^{-\ga}x$ and
  \[
  A_p\ph(x) \le -Kx + Cx^\al + Cx^{-1} - Kx^{\al-2}.
  \]
If $x<p^\ga\ell/(1-c_2)$, then $g(p^{-\ga}x)=p^{-\ga}x-\ell$ and
  \[
  A_p\ph(x) \le Cx^\al + \frac{x - p^\ga\ell}{xp^\ga\ell}
    - Kx^{\al-2}
    \le Cx^\al + (p^\ga\ell)^{-1} - Kx^{\al-2}.
  \]
But in this case, $(p^\ga\ell)^{-1}<Cx^{-1}$. It therefore follows
that
  \[
  A_p\ph(x) \le C - Kx - Kx^{\al-2}
  \]
for all $x\in E_p$.

Let $\ep>0$. Define
  \[
  L = \sup_{p<1/2} \,\, \sup_{x\in E_p} A_p\ph(x) < \infty
  \]
and let $m=L(1-\ep)/\ep$. Choose $M>0$ such that
$x\notin[M^{-1},M]$ implies $A_p\ph(x)<-m$ for all $p<1/2$. By
Corollary 4.9.8 in \citep{EK},
  \[
  \eta_p([M^{-1},M]) \ge \eta_p(\{x : A_p\ph(x) \ge -m\})
    \ge \frac m {L + m} = 1 - \ep.
  \]
The family of measures $\{\eta_p\}$ is therefore relatively
compact on $(0,\infty)$. By passing to a subsequence, we can
assume that $\eta_p\To\eta$ for some probability measure $\eta$ on
$(0,\infty)$.

Now let $p^\ga W_0$ have distribution $\eta_p$ and let $Z_p$ be
given by \eqref{Zp}. By Theorem \ref{T:Poisson}, $Z_p\To Z$, where
$Z$ satisfies \eqref{Z} with $PZ(0)^{-1}=\eta$. Fix
$t_1\le\cdots\le t_n$. Then
  \begin{align*}
  (Z_p(t_1), \ldots, Z_p(t_n))
    &= p^\ga(W_{\flr{t_1p^{-1}}},
    \ldots, W_{\flr{t_np^{-1}}})\\
  &\eqd p^\ga(W_0,
    W_{\flr{t_2p^{-1}} - \flr{t_1p^{-1}}},
    \ldots,
    W_{\flr{t_np^{-1}} - \flr{t_1p^{-1}}})\\
  &= (Z_p(0), Z_p(t_2 - t_1), \ldots, Z_p(t_n - t_1)) + \ep,
  \end{align*}
where $\ep_j=Z_p(h_j)-Z_p(t_j-t_1)$ and
$h_j=(\flr{t_jp^{-1}}-\flr{t_1p^{-1}})p$. Note that $h_j\to
t_j-t_1$ as $p\to0$ and, for fixed $t$, $Z$ is almost surely
continuous at $t$. Hence, $\ep\to0$ a.s., which gives
  \[
  (Z_p(t_1), \ldots, Z_p(t_n))
    \To (Z(0), Z(t_2 - t_1), \ldots, Z(t_n - t_1)).
  \]
But
  \[
  (Z_p(t_1), \ldots, Z_p(t_n))
    \To (Z(t_1), \ldots, Z(t_n)),
  \]
so $Z$ is a stationary process, and $\eta$ is a stationary
distribution for $Z$. The uniqueness of $\eta$ follows from Lemma
\ref{L:Zprop}. \qed

For the proof of Theorem \ref{T:be<1}, note that $\xih_p$ is a
continuous time Markov chain on
$E_p=[p^{-\tau}(p^\ga\ell-c_p),\infty)$ with generator
  \begin{equation}\label{gen}
  \begin{split}
  A_p\ph(x) &= p^{-\nu+1}
    (\ph(x - p^{\ga-\tau}g(p^{\tau-\ga}x + p^{-\ga}c_p))
    - \ph(x))\\
  &\quad + p^{-\nu}(1 - p)
    (\ph(x + p^{\ga-\tau}c_1(p^{\tau-\ga}x + p^{-\ga}c_p)^\al)
    - \ph(x)).
  \end{split}
  \end{equation}
We will use the same argument as in the proof of Theorem
\ref{T:be=1}, this time using the Lyapunov function
$\ph(x)=|x|^r$, where $r$ is sufficiently large. Our key estimate
on $A_p\ph(x)$ is given in the following lemma and is valid as
long as $|x|$ is not too large.

\begin{lemma}\label{L:genest1}
Suppose $\be<1$. Let $\ph(x)=|x|^r$, where $r\ge 2$, and let $A_p$
be given by \eqref{gen}. Let $0<\de<M<\infty$ be arbitrary.  Then
there exists $p_0>0$ and strictly positive, finite constants $C$
and $K$ such that
  \[
  A_p\ph(x) \le C - K|x|^r
  \]
for all $p\le p_0$ and all $x\in E_p$ satisfying $\de\le p^\tau
x+c_p\le M$.
\end{lemma}

\pf For notational simplicity, let us define $y_p(x)=p^\tau x+c_p$
so that
  \begin{equation*}
  \begin{split}
  A_p\ph(x) &= p^{-\nu+1}
    (\ph(x - p^{\ga-\tau}g(p^{-\ga}y_p)) - \ph(x))\\
  &\quad + p^{-\nu}(1 - p)
    (\ph(x + p^{\ga-\tau}c_1(p^{-\ga}y_p)^\al) - \ph(x)).
  \end{split}
  \end{equation*}
Either $g(x)=c_2x^\be$ or $g(x)<c_2x^\be$. Note that there exists
$x_0>\ell$ such that $g(x)=c_2x^\be$ if and only if $x\ge x_0$.
Hence, if $g(p^{-\ga}y_p)<c_2(p^{-\ga}y_p)^\be$, then
$p^{-\ga}y_p<x_0$, which implies $x<p^{-\tau}(p^\ga x_0-c_p)$. If
$p$ is sufficiently small, this implies $x<0$. Since $\ph$ is
decreasing on $(-\infty,0]$, it follows that
  \begin{equation*}
  \begin{split}
  A_p\ph(x) &\le p^{-\nu+1}(\ph(x - p^{\ga-\tau-\ga\be}c_2y_p^\be)
    - \ph(x))\\
  &\quad + p^{-\nu}(1 - p)(\ph(x + p^{\ga-\tau-\ga\al}c_1y_p^\al)
    - \ph(x))
  \end{split}
  \end{equation*}
for all $x\in E_p$.

Observe that
  \begin{align*}
  |\ph(z) - \ph(x) - \ph'(x)(z - x)|
    &= \left|{\int_x^z (z - u)\ph''(u)\,du}\right|\\
  &\le C|z - x|^2(|x|^{r-2} + |z|^{r-2})\\
  &\le C|x|^{r-2}|z - x|^2 + C|z - x|^r.
  \end{align*}
Hence,
  \begin{equation*}
  \begin{split}
  A_p\ph(x) &\le -\ph'(x)p^{-\tau}(p^{-\nu+1+\ga-\ga\be}c_2y_p^\be
    - p^{-\nu+\ga-\ga\al}c_1(1 - p)y_p^\al)\\
  &\quad + C|x|^{r-2}(p^{-\nu+1+2\ga-2\tau-2\ga\be}c_2^2y_p^{2\be}
    + p^{-\nu+2\ga-2\tau-2\ga\al}c_1^2y_p^{2\al})\\
  &\quad + C(p^{-\nu+1+r\ga-r\tau-r\ga\be}c_2^ry_p^{r\be}
    + p^{-\nu+r\ga-r\tau-r\ga\al}c_1^ry_p^{r\al}).
  \end{split}
  \end{equation*}
We can simplify these exponents by observing that
  \begin{align*}
  -\nu + \ga - \ga\al &= 0\\
  -\nu + 1 + \ga - \ga\be &= 0\\
  -\nu + 2\ga - 2\tau - 2\ga\al &= 1\\
  -\nu + 1 + 2\ga - 2\tau - 2\ga\be &= 0\\
  -\nu + 1 + r\ga - r\tau - r\ga\be &= \tau(r - 2)\\
  -\nu + r\ga - r\tau - r\ga\al &= r - 1 + \tau(r - 2).
  \end{align*}
Thus,
  \begin{equation*}
  \begin{split}
  A_p\ph(x)
    &\le -\ph'(x)p^{-\tau}(c_2y_p^\be - c_1(1 - p)y_p^\al)
    + C|x|^{r-2}(y_p^{2\be} + py_p^{2\al})\\
  &\quad + C(p^{\tau(r-2)}y_p^{r\be}
    + p^{r-1+\tau(r-2)}y_p^{r\al}).
  \end{split}
  \end{equation*}
Since $\ph'(x)$ and $c_2y_p^\be-c_1(1-p)y_p^\al$ have the same
sign, this gives
  \begin{equation}\label{genest}
  \begin{split}
  A_p\ph(x)
    &\le -r|x|^{r-1}p^{-\tau}|c_2y_p^\be - c_1(1 - p)y_p^\al|
    + C|x|^{r-2}(y_p^{2\be} + py_p^{2\al})\\
  &\quad + C(p^{\tau(r-2)}y_p^{r\be}
    + p^{r-1+\tau(r-2)}y_p^{r\al})
  \end{split}
  \end{equation}
for all $x\in E_p$.

If $r\ge2$ and $\de\le y_p\le M$, then
  \[
  A_p\ph(x) \le -r|x|^{r-1}p^{-\tau}c_2y_p^\al
    |y_p^{\be-\al} - c_p^{\be-\al}| + C|x|^{r-2} + C.
  \]
By the Mean Value Theorem,
  \begin{align*}
  \psi_p(x)
    &\le -K|x|^{r-1}p^{-\tau}|y_p - c_p| + C|x|^{r-2} + C\\
  &= -K|x|^r + C|x|^{r-2} + C,
  \end{align*}
which completes the proof. \qed

The following two lemmas provide the needed estimates on $A_p\ph$
in the extreme regimes.

\begin{lemma}\label{L:genest2}
Suppose $\be<1$. Let $\ph(x)=|x|^r$, where $r\ge 2$, and let $A_p$
be given by \eqref{gen}. Then there exists $p_0>0$, $M<\infty$ and
$K>0$ such that
  \[
  A_p\ph(x) \le - K|x|^{(r-1)\wedge(r-1+\be)}
  \]
for all $p\le p_0$ and all $x\in E_p$ satisfying $p^\tau x+c_p>M$.
\end{lemma}

\pf Let $p\le p_0$ and $y_p=p^\tau x+c_p>M$. If $p_0$ is
sufficiently small and $M$ is sufficiently large, then $x\ge
Kp^{-\tau}$ and $y_p\le x$. By \eqref{genest},
  \begin{align*}
  A_p\ph(x) &\le -K|x|^{r-1}y_p^\be + C|x|^{r-2}y_p^{2\be}
    + Cy_p^{r\be}\\
  &= -|x|^{r-1}y_p^\be(K - C|x|^{-1}y_p^\be
    - C|x|^{-r+1}y_p^{\be(r-1)}).
  \end{align*}
If $\be\le0$, then for $p$ sufficiently small,
  \[
  A_p\ph(x) \le -|x|^{r-1}y_p^\be(K - C|x|^{-1} - C|x|^{-r+1})
    \le -K|x|^{r-1+\be}.
  \]
If $\be>0$, then
  \[
  A_p\ph(x) \le -|x|^{r-1}y_p^\be(K - C|x|^{\be-1}
    - C|x|^{(\be-1)(r-1)}),
  \]
so for $p$ sufficiently small,
$A_p\ph(x)\le-K|x|^{r-1}y_p^\be\le-K|x|^{r-1}$. \qed

\begin{lemma}\label{L:genest3}
Suppose $\be<1$. Let $\ph(x)=|x|^r$, where $r\ge 2$, and let $A_p$
be given by \eqref{gen}. Then there exists $p_0>0$, $\de>0$ and
$K>0$ such that
  \[
  A_p\ph(x) \le -K|x|^{r\wedge(r-2\al/(1-\be))}
  \]
for all $p\le p_0$ and all $x\in E_p$ satisfying $p^\tau
x+c_p<\de$.
\end{lemma}

\pf Let $p\le p_0$ and $y_p=p^\tau x+c_p<\de$. Note that since
$x\in E_p$, $y_p\ge p^\ga\ell$. If $p_0$ and $\de$ are
sufficiently small, then $x<0$ and $Kp^{-\tau}\le|x|\le
Cp^{-\tau}$. By \eqref{genest}, for $\de$ sufficiently small,
  \begin{equation*}
  \begin{split}
  A_p\ph(x) &\le -|x|^ry_p^\al(K|y_p^{\be-\al} - c_p^{\be-\al}|
    - C(p^{2\tau}y_p^{2\be-\al}
    + p^{2\tau+1}y_p^\al)\\
  &\quad - C(p^{\tau r+\tau(r-2)}y_p^{r\be-\al}
    + p^{\tau r+r-1+\tau(r-2)}y_p^{r\al-\al}))\\
  &\le -|x|^ry_p^\al(K - C(p^{2\tau}y_p^{2\be-\al}
    + p^{2\tau(r-1)}y_p^{r\be-\al})\\
  &\quad - C(p^{2\tau+1}y_p^\al
    + p^{(2\tau+1)(r-1)}y_p^{\al(r-1)})).
  \end{split}
  \end{equation*}
Let us first estimate the term $p^{2\tau}y_p^{2\be-\al}$. If
$2\be-\al\ge0$, then $p^{2\tau}y_p^{2\be-\al}\le Cp^{2\tau}$. If
$2\be-\al<0$, then $p^{2\tau}y_p^{2\be-\al}\le
Cp^{2\tau+\ga(2\be-\al)}$. Note that $2\tau+\ga(2\be-\al)=\ga+1$.
Hence, for all values of $\al$ and $\be$, there exists some $s>0$
such that $p^{2\tau}y_p^{2\be-\al}\le p^s$.

Similarly, for the remaining terms in the above inequality, we
observe that
  \begin{align*}
  2\tau(r-1) + \ga(r\be - \al) &= (2\tau + \ga\be)(r - 1) + 1
    = \ga(r - 1) + 1\\
  2\tau + 1 + \ga\al &= \ga\\
  (2\tau + 1)(r - 1) + \ga\al(r - 1) &= \ga(r - 1).
  \end{align*}
Therefore, if $p_0$ is sufficiently small, then
$A_p\ph(x)\le-K|x|^ry_p^\al$. If $\al<0$, then
$A_p\ph(x)\le-K|x|^r$. If $\al\ge0$, then
  \[
  A_p\ph(x) \le -K|x|^rp^{\ga\al} \le -K|x|^{r-\ga\al/\tau}.
  \]
Since $\ga\al/\tau=2\al/(1-\be)$, this completes the proof. \qed

\bigskip

\noindent{\bf Proof of Theorem \ref{T:be<1}.} Suppose $\be<1$ and
$\eta_p$ is the stationary distribution for $\{p^{-\tau}(p^\ga
W_n-c_p)\}$. Then $\eta_p$ is the stationary distribution for
$\xih_p$. Let $\ph(x)=|x|^r$, where $r\ge2$. By Lemmas
\ref{L:genest1}, \ref{L:genest2}, and \ref{L:genest3}, if $r$ is
sufficiently large, there exists $p_0>0$ and strictly positive,
finite constants $C$ and $K$ such that
  \[
  A_p\ph(x) \le C - K|x|^s
  \]
for some $s>0$ and all $p\le p_0$ and $x\in E_p$. As in the proof
of Theorem \ref{T:be=1}, this implies that the family of measures
$\{\eta_p\}$ is relatively compact on $\RR$. By passing to a
subsequence, we can assume that $\eta_p\To\eta$ for some
probability measure $\eta$ on $\RR$.

Let $p^{-\tau}(p^\ga W_0-c_p)$ have distribution $\eta_p$, let
$Z_p$ be given by \eqref{Zp}, and let $\xi_p$ be given by
\eqref{xip} with $\ze_p\equiv c_p$. Note that $\xi_p(0)$ converges
in distribution, so $p^\tau\xi_p(0)=Z_p(0)-\ze_p(0)\to0$ in
probability. Hence, by Theorem \ref{T:CLT}, $\xi_p\To\xi$, where
$\xi$ satisfies \eqref{OU} with $P\xi(0)^{-1}=\eta$. As in the
proof of Theorem \ref{T:be=1}, $\xi$ is a stationary process, so
$\eta$ is the stationary distribution for $\xi$. \qed

\section*{Acknowledgements}

The authors gratefully acknowledge Tom Kurtz, Ruth Williams, and
Timo Sepp\"al\"ainen for their support and helpful ideas.

\nocite{*}
\bibliographystyle{plain}
\bibliography{genTCP2}

\end{document}